\documentclass[11pt, amsthm]{elsart}

  \usepackage{ccaption}

  \captionnamefont{\bfseries}
  \captiontitlefont{\small\sffamily}
  \captiondelim{ --- }
  \hangcaption
  \renewcommand{\figurename}{Table}

\usepackage{amssymb}
\usepackage{amsmath}
\usepackage{enumerate}
\usepackage{mathrsfs}
\usepackage{graphics}
\usepackage{cases}
\usepackage{bbding}
\usepackage{verbatim}
\usepackage{pifont}

\theoremstyle{plain}

\newtheorem{theorem}{Theorem}[section] 
\newtheorem{lemma}[theorem]{Lemma}     
\newtheorem{corollary}[theorem]{Corollary}
\newtheorem{proposition}[theorem]{Proposition}

\theoremstyle{definition}

\newtheorem{definition}[theorem]{Definition}

\newtheorem{remark}[theorem]{Remark}
\newtheorem{example}[theorem]{Example}

\newcommand{\vspc}{\vspace{2mm}}

\newcommand{\im}{\mbox{\rm im}}

\newcommand{\gl}{\mathcal{L}}
\newcommand{\gr}{\mathcal{R}}

\newcommand{\gh}{\mathcal{H}}

\begin{document}

\begin{frontmatter}

\title{ On Residual Finiteness of Monoids,  their Sch\"{u}tzenberger Groups and \\ Associated Actions}

\author{R. Gray}\footnote{\textit{Address:} 
School of Mathematics, University of East Anglia, Norwich NR4 7TJ, U.K.,
\textit{e-mail:} \textup{\texttt{Robert.D.Gray@uea.ac.uk.}} 
}
\and
\author{N. Ru\v{s}kuc}\footnote{\textit{Address:} School of Mathematics and Statistics, University of St Andrews, St Andrews KY16 9SS, U.K., \textit{e-mail:} \textup{\texttt{nik@mcs.st-and.ac.uk.}}}

\maketitle

\begin{abstract}
%
%
%
%
%
%
In this paper we discuss connections between
the following properties:
(RFM) residual finiteness of a monoid $M$; (RFSG) residual finiteness of Sch\"{u}tzenberger groups
of $M$;  and (RFRL) residual finiteness of the natural actions of $M$ on its Green's $\mathcal{R}$- and $\mathcal{L}$-classes.
The general question is whether (RFM) implies (RFSG) and/or (RFRL), and vice versa.
We consider these questions in all the possible combinations of the following situations:
$M$ is an arbitrary monoid; $M$ is an arbitrary regular monoid;
every $\mathcal{J}$-class of $M$ has finitely many $\mathcal{R}$- and $\mathcal{L}$-classes;
$M$ has finitely many left and right ideals.
In each case we obtain complete answers, which are summarised in Table~\ref{table1}.  

\bigskip


\end{abstract}


\end{frontmatter}

\section{Introduction}
\label{sec1}
Residual finiteness is one of the most important algebraic finiteness conditions, and has been widely studied in the context of groups, semigroups, monoids and other algebraic structures.
In group theory, following the famous results of 
Malcev \cite{Malcev1} and other pioneers 
such as Hirsch \cite{Hirsch1}, Gruenberg \cite{Gruenberg1} and Hall \cite{PHall1},
there is now a large body of literature devoted to this property;
for an early survey see \cite{Magnus1}.
In semigroup theory early work had a similar nature as in groups;
for example see \cite{Evans}, \cite{Carlisle}, \cite{Lallement}, \cite{Golubov3}, \cite{Golubov4}, \cite{Golubov2}, \cite{Golubov1}, \cite{GolubovSapir}.
Samples of more recent work can be found in \cite{Auinger}, \cite{LallementRosaz}, \cite{Gray1}, \cite{Ruskuc&Thomas}. The work on residual finiteness in other algebraic structures includes rings and modules (\cite{Varadarajan}, \cite{Faith}), Lie algebras (\cite{PremetSemenov}, \cite{Zaicev}, \cite{BahturinZaicev}), lattices (\cite{AdamsGratzer}), and universal algebra (\cite{Evans1}, \cite{BanaschewskiNelson}, \cite{GolubovSapir1}, \cite{Gray2}).

The purpose of this paper is to discuss residual finiteness for monoids.
Associated with an arbitrary monoid $M$ are a number of natural
actions. A monoid acts on its Green's $\gl$-classes via right
multiplication, and, dually, 
on its $\gr$-classes via left multiplication  (see Section \ref{sec2} for background on Green's relations). Of course $M$ also acts on its
elements by right multiplication. Arising from this action is a family
of groups, called the Sch\"{u}tzenberger groups of $M$, found by taking the
action of the set-wise stabiliser of an $\gh$-class $H$, restricted to $H
$ and made faithful by factoring out the kernel.
These groups encode much of the group theoretic
structure of $M$. For instance, every maximal subgroup of $M$ arises as
a Sch\"{u}tzenberger group. The obvious connections between actions and
congruences leads one to consider the relationship between residual
finiteness of $M$ and the properties of the natural actions associated
with $M$ outlined above. More precisely, here we will explore the
connections between a monoid $M$ being residually finite and the
residual finiteness of:
\begin{itemize}
\item[--]
the Sch\"{u}tzenberger groups of $\mathcal{H}$-classes of $M$;
\item[--]
the action of $M$ on its $\mathcal{L}$-classes;
\item[--]
the (left) action of $M$ on its $\mathcal{R}$-classes.
\end{itemize}
\sloppypar{Again, we refer the reader to Section \ref{sec2} for an introduction to action, Sch\"{u}tzenberger groups and residual finiteness. }

For a while during this investigation the authors had hoped that the following general assertion might be true:
A monoid is residually finite if and only if its Sch\"{u}tzenberger groups, and its actions on $\mathcal{R}$- and $\mathcal{L}$-classes are all residually finite.
For instance, Golubov  \cite{Golubov1} proves that a \emph{regular}
monoid in which each $\mathcal{J}$-class $J$ has only finitely many $\mathcal{R}$- and $\mathcal{L}$-classes is residually finite if and only if all its maximal subgroups are residually finite.
Here the maximal subgroups \emph{are} the Sch\"{u}tzenberger groups, while the finiteness condition on $\mathcal{J}$-classes (we will call this \emph{finite type}) implies residual finiteness of the actions on $\mathcal{R}$- and $\mathcal{L}$-classes (see Theorem \ref{thm10}).
Sadly, our hopes were premature: the desired general result does not hold in either direction, as we will see below.

In this paper we will systematically compare and contrast regular and non-regular monoids, in the following situations:
\begin{itemize}
\item[--]
general monoids (Section \ref{sec3});
\item[--]
monoids in which every $\mathcal{J}$-class is of finite type (Sections \ref{sec4} and \ref{sec4b});
\item[--]
monoids with finitely many left- and right ideals (Section \ref{sec5}).
\end{itemize}

Section \ref{sec4a} contains a couple of further auxiliary results, concerning ideals and Rees quotients.
We state our results for monoids, for the sake of convenience of having an identity element around.
The results remain valid verbatim for semigroups, which can be seen by utilising the standard device of adjoining an identity element to a semigroup.

The results of the paper are summarised in Table~\ref{table1}.
\begin{figure}
\begin{center}
\begin{tabular}{|l|l|c|c|c|} \hline
\multicolumn{2}{|c|}{Properties of $M$}& 
\parbox{25mm}{\centering $M$ r.f.$\Rightarrow$ \\ all $\Gamma(H)$ r.f.} & 
\parbox{25mm}{\centering $M$ r.f. $\Rightarrow$\\ $M/\mathcal{L}$ r.f.} &
\parbox{30mm}{\ \\  \centering all $\Gamma(H)$ r.f. \&\\ $M/\mathcal{L}, M/\mathcal{R}$ r.f.\\ $\Rightarrow$ $M$ r.f.\vspace{1mm}}
\\
\hline\hline
general &  non-reg & 
\parbox{25mm}{\ \\ \centering \checkmark\\ (Thm \ref{thm1})\vspace{1mm}} & 
\parbox{25mm}{\centering \ding{55} \\ (Cor \ref{corol6})} & 
\parbox{30mm}{\centering \ding{55} \\ (Cor \ref{corol9})}
\\
\cline{2-5}
&regular& 
\parbox{25mm}{\ \\  \centering \checkmark\\ (as above)\vspace{1mm}} &
\parbox{25mm}{\centering \checkmark \\ (Thm \ref{thm2})} & 
\parbox{30mm}{\centering \ding{55}\\ (Cor \ref{corol9})}
\\
\hline
\parbox{15mm}{\ \\ finite\\ $\mathcal{J}$-type} &  non-reg & 
\parbox{25mm}{\ \\  \centering \checkmark\\ (as above)\vspace{1mm}} & 
\parbox{25mm}{\centering \ding{55} \\ (Prop \ref{prop4b16})} & 
\parbox{30mm}{\centering \ding{55}\\ (Ex \ref{ex13})}
\\
\cline{2-5}
&regular& 
\parbox{25mm}{\centering \checkmark\\ (as above)} &
\parbox{25mm}{\ \\  \centering \checkmark \\ (as above \&\\Thm \ref{thm10})\vspace{1mm}} & 
\parbox{30mm}{\centering \checkmark \\ (Thm \ref{thm10})}
\\
\hline
\parbox{15mm}{\ \\finitely\\ many\\ right/left\\ ideals} &  non-reg & 
\parbox{25mm}{\ \\  \centering \checkmark\\ (as above)} & 
\parbox{25mm}{\centering \checkmark \\ (obvious)} & 
\parbox{30mm}{\centering \checkmark\\ (Thm \ref{thm16})}
\\
\cline{2-5}
&regular& 
\parbox{25mm}{\ \\  \centering \checkmark\\ (as above)\vspace{1mm}} &
\parbox{25mm}{\centering \checkmark \\ (as above)} & 
\parbox{30mm}{\centering \checkmark \\ (as above)}
\\
\hline
\end{tabular} 
\end{center}
\caption{Summary of results.}
\label{table1}
\end{figure}
\section{Preliminaries}
\label{sec2}
A monoid $M$ is said to be \emph{residually finite} if for any two distinct elements $x,y\in S$ there exists a finite monoid $N$ and a homomorphism $f\::\: M\rightarrow N$ such that $f(x)\neq f(y)$. 
This is equivalent to the following properties:
\begin{itemize}
\item[--]
there exists a congruence $\rho$ of finite index (i.e. with finitely many classes) on $M$ such that
 $x/\rho\neq y/\rho$;
\item[--]
the intersection of all congruences of finite index on $M$ is trivial.
\end{itemize}

Green's equivalences $\mathcal{R}$, $\mathcal{L}$,  $\mathcal{J}$, $\mathcal{H}$ and $\mathcal{D}$ are defined as follows:
\begin{eqnarray*}
& \mathcal{R}=\{ (x,y)\::\: xM=yM\},\ \mathcal{L}=\{(x,y)\::\: Mx=My\},\\
&\mathcal{J}=\{(x,y)\::\: MxM=MyM\},\ 
\mathcal{H}=\mathcal{R}\cap\mathcal{L},\ \mathcal{D}=\mathcal{R}\circ\mathcal{L}=\mathcal{L}\circ\mathcal{R}.&
\end{eqnarray*}
Also for a set $X$ we let
\[
\Delta_X= \{ (x,x)\::\: x\in X\},\ \Phi_X=\{ (x,y)\::\: x,y\in X\}.
\]
A monoid $M$ is said to be \emph{regular} if for every $x\in M$ there exists $y\in M$ such that $xyx=x$.
In a regular monoid every $\mathcal{R}$-class and every $\mathcal{L}$-class contain at least one idempotent.
If an $\mathcal{H}$-class contains an idempotent, then it is a group; these are precisely the maximal subgroups of $M$.

Clearly, $\mathcal{R}$-classes (respectively, $\mathcal{L}$- and $\mathcal{J}$-classes) are in one-one correspondence with principal right (resp. left, two-sided) ideals.
Thus, the condition that the monoid $M$ should have finitely many left- and right ideals, which we will consider in Section \ref{sec5}, is equivalent to there being finitely many $\mathcal{L}$- and $\mathcal{R}$-classes (and hence also finitely many $\mathcal{J}$-classes). A weaker condition, which will be under consideration in Section \ref{sec4} is that every $\mathcal{J}$-class has \emph{finite type}, by which we mean that it contains only finitely many $\mathcal{R}$- and $\mathcal{L}$-classes.

A \emph{(right) action} of a monoid $M$ on a set $X$ is a mapping $X\times M\rightarrow X$, $(x,m)\mapsto xm$, such that $(xm)n=x(mn)$ and $x1=x$ for all $x\in X$, $m,n\in M$.
In element $x_0\in X$ such that $x_0M=X$ is called a \emph{source} for the action.
If $\rho$ is a right congruence on $M$, then $M$ acts on the set $M/\rho$ of all $\rho$-classes by
$(x/\rho)m=(xm)/\rho$, and $1/\rho$ is a source.
Conversely, if $M$ acts on $X$ with source $x_0$, the relation $\sim$ defined by $\{ (m,n)\::\: x_0m=x_0n\}$ is a right congruence; we say that $\sim$ is \emph{determined} by the action. 
The above two constructions are mutually inverse.
There is also a natural two-sided congruence $\{ (m,n)\::\: (\forall x\in X)(xm=xn)\}$ 
associated with the action; we call it the \emph{kernel} of the action.

Suppose $M$ acts on two sets $X$ and $Y$. A mapping $f\::\: X\rightarrow Y$ is a \emph{homomorphism of actions} if $f(xm)=f(x)m$ for all $x\in X$, $m\in M$.
An action of a monoid $M$ on a set $X$ is \emph{residually finite} if for any two distinct $x,y\in X$ there exists an action of $M$ on a finite set $Y$ and a homomorphism $f\::\: X\rightarrow Y$ such that
$f(x)\neq f(y)$.
It is a routine exercise to show that for actions with a source point this is equivalent to the condition that the right congruence determined by the action is the intersection of finite index right congruences. All the above notions have duals for left actions and congruences.

Green's equivalence $\mathcal{L}$ (respectively $\mathcal{R}$) is a right (resp. left) congruence,
and so there is a natural action (resp. left action) of $M$ on its $\mathcal{L}$-
(resp. $\mathcal{R}$-) classes.

It is possible to associate a group to an \emph{arbitrary} $\gh$-class $H$ in a monoid $M$, regardless of whether it contains an idempotent.  
First define the \emph{stabilizer} of $H$ in $M$ under the right multiplication action of $M$ on itself:
$\mathrm{Stab}(H) = \{ s \in S : Hs=H \}$.
Clearly, $\mathrm{Stab}(H)$ is a submonoid of $M$, and it acts by right multiplication on $H$.
The kernel $\sigma=\sigma(H)$  of this action is given by: 
$(x,y) \in \sigma$ if and only if $hx = hy$
for all $h\in H$. 
It turns out that $\Gamma(H)=\mathrm{Stab}(H) / \sigma$ is a group, 
called the \emph{Sch\"{u}tzenberger group} of $H$.
The group $\Gamma(H)$ acts naturally on $H$ via $h\cdot (s/\sigma)=hs$.
The standard properties of Sch\"{u}tzenberger groups are summarized in the following result (see \cite[Section~2.3]{Lallement2} for proofs of these facts). 

\begin{proposition}
\label{schutzstabproperties}
Let $M$ be a monoid, let $H$ be an
$\gh$-class of $M$, and let 
$h \in H$ be an arbitrary element. Then:
\begin{enumerate}[(i)]
\item
\label{item-sch1}
$\mathrm{Stab}(H) = \{ s \in M : hs \in H  \}$.
\item 
\label{item-sch2}
$\sigma(H) = \{ (u,v) \in \mathrm{Stab}(H) \times \mathrm{Stab}(H) : hu = hv  \}$ (i.e. the kernel of the action of $\mathrm{Stab}(H)$ on $H$ is equal to the right congruence associated with the action).
\item 
\label{item-sch3}
$H = h \mathrm{Stab}(H)=h\cdot \Gamma(H)$.
\item 
\label{item-sch4}
If $H$ and $H'$ belong to the same $\gl$-class of $M$ then $\mathrm{Stab}(H) = \mathrm{Stab}(H')$. 
\item 
\label{item-sch5}
The action of $\Gamma(H)$ on $H$ is regular;
in particular $|\Gamma(H)| = |H|$.
\item 
\label{item-sch6}
If $H_1$ is an $\gh$-class of $M$ belonging to the same $\gr$-class
(or to the same $\gl$-class) as $H$ then $\Gamma(H_1) \cong \Gamma(H)$.
\item 
\label{item-sch7}
If $H$ is a group then $\Gamma(H) \cong H$. 
\end{enumerate}
\end{proposition}

Of course, one could left-right dualise the definition of the Sch\"{u}tzenberger group,
to obtain its left version $\Gamma_l(H)$; however, it turns out that $\Gamma_l(H)\cong\Gamma(H)$.

\section{General Monoids}
\label{sec3}

In this section we ask: Does residual finiteness of a monoid imply residual finiteness of its Sch\"{u}tzenberger groups and of its actions on $\mathcal{R}$- and $\mathcal{L}$-classes? 
Do the last two conditions imply the first?

It is clear that residual finiteness is a hereditary property: 
every substructure of a residually finite structure is residually finite. In particular, maximal subgroups of residually finite monoids are residually finite.
This can be generalised to Sch\"{u}tzenberger groups:

\begin{theorem}[Gray, Ru\v{s}kuc \cite{Gray1}]
\label{thm1}
All Sch\"{u}tzenberger groups of a residually finite monoid are residually finite.
\end{theorem}

\begin{proof}
We provide a very brief sketch here, for the sake of completeness:
Take two distinct elements $x/\sigma,y/\sigma\in \Gamma(H)$. Then $hx\neq hy$ in $M$,
where $h\in H$ is arbitrary.
Let $M/\rho$ be a finite quotient of $M$ in which $(hx)/\rho\neq (hy)/\rho$.
Let $\overline{H}$ be the $\gh$-class of $h/\rho$ in $M/\rho$.
The obvious homomorphism $\Gamma(H)\rightarrow\Gamma(\overline{H})$ separates $x/\sigma$ and $y/\sigma$.
\end{proof}

Now we turn our attention to the actions on $\mathcal{R}$- and $\mathcal{L}$-classes. For regular monoids we have:

\begin{theorem}
\label{thm2}
If $M$ is a regular, residually finite monoid then the actions of $M$ on its $\mathcal{R}$- and $\mathcal{L}$-classes are also residually finite.
\end{theorem}

\begin{proof}
We will prove the assertion for $\mathcal{L}$-classes; the proof for $\mathcal{R}$-classes is dual.
Let $s,t\in M$ be any elements with $(s,t)\not\in \mathcal{L}$.
We need to find a finite index right congruence $\rho$ such that $\mathcal{L}\subseteq \rho$
and $(s,t)\not\in\rho$.
Since $M$ is regular it follows that there exist idempotents $e,f\in M$ such that $s\mathcal{L}e$ and $t\mathcal{L} f$.
It is known that two idempotents $g$ and $h$ are $\mathcal{L}$-related if and only if $gh=g$ and $hg=h$
(\cite[Exercise 2.6.3]{Howie1}); in particular, $ef\neq e$.
Since $M$ is residually finite, there exists a homomorphism $\phi : M\rightarrow N$, $N$ finite,
with $\phi(ef)\neq \phi(e)$.
Let $\rho$ be the pre-image of the $\mathcal{L}$ equivalence on $N$, which we will denote by $\mathcal{L}^N$:
$$
\rho=\phi^{-1}(\mathcal{L}^N)=\{ (x,y)\in M\times M\::\: (\phi(x),\phi(y))\in\mathcal{L}^N\}.
$$
Clearly, since $x\mathcal{L}y$ implies $\phi(x)\mathcal{L}^N\phi(y)$, we have $\mathcal{L}\subseteq \rho$.
Finiteness of $N$ implies that $\rho$ has finite index.
The elements $\phi(e)$, $\phi(f)$ are idempotents, and from $\phi(e)\phi(f)\neq \phi(e)$ it follows that $(\phi(e),\phi(f))\not\in \mathcal{L}^N$.
Thus $(e,f)\not\in\rho$, and the proof is complete.
\end{proof}

\begin{remark}
The above argument actually 
shows that any regular residually finite monoid has the following property:
for any two elements $s,t\in M$ with $(s,t)\not\in \mathcal{L}$ there exists a homomorphism
$f:M\rightarrow N$, $N$ finite, such that $f(s)$ and $f(t)$ are not $\mathcal{L}$-related in $N$.
Let us call this property \emph{$\mathcal{L}$-separation}.
One could think (as the authors did for a while) that $\mathcal{L}$-separation is equivalent to the action of
$M$ on its $\mathcal{L}$-classes being residually finite.
This, however, is not the case as the following example shows.
\end{remark}

\begin{example}
Define a monoid $M$ with zero as follows:
$$
M=\{ a^i\::\: i\geq 0\} \cup \{ b_j\::\: j\in\mathbb{Z}\}\cup\{0\};
$$
$\{ a^i\::\: i\geq 1\}$ form a free semigroup of rank $1$; 
the remaining multiplications are governed by
$$
a^ib_j=b_{i+j},\ b_ja^i=b_jb_k=0.
$$
Note that $M$ is generated by $a$ and $b_0$.
We claim that $M$ has a residually finite action on $\gl$-classes, but
does not have the $\gl$-separation property. Moreover, we claim that $M$ itself is also
residually finite.

To show that $M$ is residually finite,
consider the monoids
$$
T_n=\langle a,b\:|\: a^{n+1}=a,\ a^nb=b,\ ba=b^2=0\rangle\ (n\geq 1).
$$
Note that $T_n$ is finite of order $2n+2$, and that there is a natural homomorphism
$\theta_n\::\: M\rightarrow T_n$, $a\mapsto a$, $b_0\mapsto b$.
Now let $s,t\in S$, $s\neq t$. If $s$ and $t$ come from two different sets
$\{1\}$, $\{ a^i\::\: i\geq 1\}$, $\{ b_j\::\: j\in\mathbb{Z}\}$ or $\{0\}$, then
$\theta_1(s)\neq \theta_1(t)$.
If $s=a^i$, $t=a^j$ with $i<j$ then $\theta_j(s)\neq \theta_j(t)$.
Finally, if $s=b_k$, $t=b_l$ with $k<l$ then $\theta_{l-k+1}(s)\neq \theta_{l-k+1}(t)$.

Now we show that $b_0$ and $b_1$, which are clearly not $\mathcal{L}$-related in $M$,
cannot be $\mathcal{L}$-separated in a finite homomorphic image of $S$.
Let $T$ be a finite monoid and let $\theta\::\: M\rightarrow T$ be a homomorphism.
Let $m\in\mathbb{N}$ be such that $\theta(a^m)=\theta(a^{2m})$.
Then we have
$$
\theta(b_1)=\theta(ab_0)=\theta(a)\theta(b_0),
$$
and
\begin{eqnarray*}
\theta(b_0)=\theta(a^mb_{-m}) & = & \theta(a^m)\theta(a^mb_{-m}) \\
 & = & \theta(a^m)\theta(b_0) 
 =  \theta(a^{m-1})\theta(ab_0)=\theta(a^{m-1})\theta(b_1);
\end{eqnarray*}
hence $\theta(b_0)\mathcal{L}^T \theta(b_1)$.
\end{example}

The assertion of Theorem \ref{thm2} does not extend to non-regular monoids. In order to construct a counter-example, let us
introduce a construction. 
(An alternative example will be provided by Proposition \ref{prop4b16}.)
We start with a group $G$, and a normal
subgroup $N\unlhd G$.
We let $\overline{N}=\{ \overline{n}\::\: n\in N\}$ be a copy of $N$ disjoint from $G$.
Define a monoid $M$ by:
\begin{eqnarray*}
&M=\mathcal{M}(G,N)=\langle G,\overline{N},h\:|\: hn=\overline{n}h,\ 
he_G=e_{\overline{N}}h=h,&\\ 
&g\overline{n}=\overline{n}g=gh=h\overline{n}=0\ 
(g\in G,\ n\in N)\rangle.&
\end{eqnarray*}
The following are  easy to prove:
\begin{enumerate}[(i)]
\item
$M=\{1\}\cup G\cup \overline{N}\cup D\cup\{0\}$
where $D=\{ hg\::\: g\in G\}$,
and these are normal forms (i.e. all the above elements are distinct).
\item
The above five sets are the $\mathcal{R}$-, $\mathcal{D}$- and $\mathcal{J}$-classes of $M$.
$D$ contains no idempotents.
\item
$\mathcal{J}$ is a congruence.
\item
Two elements $hg_1$ and $hg_2$ of $D$ are $\mathcal{L}$-related if and only if $g_1$ and $g_2$ belong
to the same coset of $N$. In other words, the $\mathcal{L}$-classes of $D$ are indexed by the quotient
$G/N$.
The remaining $\mathcal{L}$-classes of $M$ are $\{1\}$, $G$, $\overline{N}$ and $\{0\}$.
\end{enumerate}

Suppose we have another normal subgroup $K\unlhd G$, and we make the quotient of $M$ induced by
factoring $G$ by $K$:
$$
Q(M,K)=\langle M\:|\: K=e_G,\ \overline{N\cap K}=e_{\overline{N}}\rangle.
$$
Then
\begin{enumerate}
\item[(v)]
$Q(M,K)\cong \mathcal{M}(G/K,NK/K)$.
\end{enumerate}
Let us prove (v). Clearly, there is a natural homomorphism $\phi$ from $Q(M,K)$ onto
$\mathcal{M}(G/K,NK/K)$, because the latter satisfies all the defining relations for the former.
This homomorphism is clearly $1-1$ on each of $\{1\}$, $G$, $\overline{N}$ and $\{0\}$ (interpreted appropriately as
constituent sets of $Q(M,K)$, rather than $M$).
Suppose that $\phi(hg_1)=\phi(hg_2)$, or, equivalently,
$h(g_1K)=h(g_2 K)$ in $\mathcal{M}(G/K,NK/K)$. Then (ii) applied to
$\mathcal{M}(G/K,NK/K)$ implies that
$g_1K=g_2K$. So 
$g_2=g_1k$ for some $k\in K$.
But then, in $Q(M,K)$ we have $hg_1=hg_1 e_G=hg_1k=hg_2$, proving (v).

\begin{proposition}
\label{prop4}
Let $G$ be a group, let $N\unlhd G$, and let $M=\mathcal{M}(G,N)$ as above.
\begin{enumerate}
\item[\textup{(i)}]
$M$ is residually finite if and only if $G$ is residually finite.
\item[\textup{(ii)}]
The action of $M$ on its $\mathcal{L}$-classes is residually finite if and only if $G/N$ is residually finite.
\end{enumerate}
\end{proposition}

\begin{proof}
(i)
($\Rightarrow$) This is immediate, since $G$ is contained in $M$.

($\Leftarrow$)
Suppose that $G$ is residually finite.
Let $s,t\in M$, $s\neq t$.
If $s$ and $t$ are not $\mathcal{J}$-related then the natural homomorphism from $M$
onto $M/\mathcal{J}$ 
separates them.
Note that the sets $I=\overline{N}\cup D\cup \{0\}$ and $J=G\cup D\cup\{0\}$ are ideals in $M$,
and that $M/I\cong \{1\}\cup G\cup\{0\}$ and $M/J\cong \{1\}\cup \overline{N}\cup\{0\}\cong \{1\}\cup N\cup\{0\}$.
Thus, if $s,t\in G$ or $s,t\in\overline{N}$ we can separate them by first mapping $M$ onto $M/I$ or $M/J$ respectively, and then onto a finite quotient using residual finiteness of $G$ (and $N$).
Finally, suppose that $s=hg_1$, $t=hg_2$, ($g_1,g_2\in G$, $g_1\neq g_2$).
Let $K$ be a finite index normal subgroup of $G$ such that $g_1K\neq g_2 K$.
But then, if $\theta$ is the natural homomorphism from $M$ to the quotient
$Q(M,K)\cong \mathcal{M}(G/K,NK/K)$, we have
$$
\theta(s)=h(g_1K)\neq h(g_2 K)=\theta(t).
$$
We conclude that $M$ is residually finite.
\vspc

(ii)
($\Rightarrow$)
We suppose that the action of $M$ on $M/\mathcal{L}$ is residually finite, and want to prove that the group $G/N$ is residually finite.
To this end, take arbitrary $g_1N,g_2N\in G/N$, $g_1N\neq g_2N$.
By (iv) we have $(hg_1,hg_2)\not\in\mathcal{L}$, and 
residual finiteness of the action implies that there exists a right congruence $\rho$ of finite index on $M$ such that $\mathcal{L}\subseteq \rho$ and $(hg_1,hg_2)\not\in\rho$.
Now $M$ acts on $M/\rho$, and so $G$ acts on $M/\rho$ as well.
In fact, since $G$ also acts on $D$, it follows that $G$ acts on the set $\{ d/\rho\::\: d\in D\}$.
Let $K$ be the kernel of this action.
Since $\rho$ has finite index, it follows that $K$ has finite index.
For $n\in N$ and $d=hg\in D$ we have
$gnN=gN$, implying $(dn,d)\in\mathcal{L}\subseteq\rho$; hence $N\leq K$.
Also, from $(hg_1,hg_2)\not\in\mathcal{L}\subseteq\rho$ it follows that $g_1K\neq g_2K$.
Hence $G/K$ is a finite quotient of $G/N$ in which $g_1N$ and $g_2N$ are separated.

($\Leftarrow$) 
Suppose $G/N$ is residually finite.
Let $s,t\in M$ with $(s,t)\not\in\mathcal{L}$.
To show residual finiteness of the action we need to find a finite index right congruence $\rho$ such that
$\mathcal{L}\subseteq \rho$ and $(s,t)\not\in\rho$.
If at least one of $s$ or $t$ does not belong to $D$ then $(s,t)\not\in\mathcal{J}$, and $\mathcal{J}$ is a finite index congruence containing $\mathcal{L}$.
So suppose that $s,t\in D$, with $s=hg_1$, $t=hg_2$.
By (iv) we have $Ng_1\neq Ng_2$.
Residual finiteness of $G/N$ and the Correspondence Theorem for groups imply that there exists
a finite index normal subgroup $K\unlhd G$ such that $N\leq K$ and $Kg_1\neq Kg_2$.
Let $\rho$ be the congruence on $M$ for which
$$
M/\rho = Q(M,K)\cong \mathcal{M}(G/K,NK/K)=\mathcal{M}(G/K,\{K\}),
$$
where $\{K\}$ stands for the trivial subgroup of $G/K$.
Finiteness of $G/K$ implies finiteness of $M/\rho$, and so $\rho$ has finite index.
From the definition of $Q(M,K)$ it is clear that $\mathcal{L}\subseteq\rho$.
Finally, $Kg_1\neq Kg_2$ and (iv) imply that $h(Kg_1)\neq h(Kg_2)$ in $\mathcal{M}(G/K,\{K\})$,
so that in the pre-image $M$ we have $(hg_1,hg_2)\not\in\rho$.
\end{proof}

\begin{corollary}
\label{corol6}
There exists a residually finite monoid which acts on its $\mathcal{L}$-classes in a non-residually-finite way.
\end{corollary}

\begin{proof}
We can take $G$ to be a non-cyclic free group (which is well known to be residually finite), take $N\unlhd G$ to be such that $G/N$ is not residually finite, form $M=\mathcal{M}(G,N)$, and apply Proposition \ref{prop4}.
\end{proof}

Are residual finiteness of the Sch\"{u}tzenberger groups and actions on $\mathcal{L}$- and $\mathcal{R}$-classes sufficient to imply residual finiteness of the monoid? Unfortunately, this is not the case, even for regular semigroups. 
In order to see this we will 
utilise the following result of Golubov concerning residual finiteness of Rees matrix semigroups over groups.
First recall the Rees matrix semigroup construction: start with a group $G$, two index sets $I$ and $J$, and a $J\times I$ matrix $P=(p_{ji})_{j\in J,i\in I}$ with entries from $G$.
The Rees matrix semigroup $S=\mathcal{M}[G;I,J;P]$ is the set $I\times G\times J$ with multiplication
$(i,g,j)(k,h,l)=(i,gp_{jk}h,l)$. For further details and information about their significance see \cite[Chapter 3]{Howie1}.
On the set $I$ define a relation $\sim_I$ by $i\sim_I k$ if and only if there exists $g\in G$ such that
$p_{ji}=p_{jk}g$ for all $j\in J$. Clearly, this is an equivalence relation; let us denote its index by $r_I$.
Analogously, on $J$ define an equivalence by $j\sim_J l$ if and only there exists $g\in G$ such that $p_{ji}=gp_{li}$ for all $i\in I$ and denote the number of its classes by $r_J$.
The \emph{rank} of $P$ is defined to be $\max(r_I,r_J)$.
If $N$ is a normal subgroup of $G$ denote by $P/N$ the matrix $(p_{ji}N)_{j\in J, i\in I}$.

\begin{theorem}[Golubov {\cite[Theorem 3]{Golubov2}}]
\label{thm8}
A Rees matrix semigroup $\mathcal{M}[G;I,J;P]$ is residually finite if and only if $G$ is residually finite and the matrix $P/N$ has finite rank for every normal subgroup $N$ of $G$ of finite index.
\end{theorem}

\begin{corollary}
\label{corol9}
There exists a regular monoid $M$ such that all its Sch\"{u}tzenberger groups (i.e. maximal subgroups) are residually finite, and so are its actions on $\mathcal{R}$- and $\mathcal{L}$-classes, but $M$ itself is not residually finite.
\end{corollary}
\begin{proof}
Let $S=\mathcal{M}[C_2;\mathbb{N},\mathbb{N};P]$, where $C_2=\{e,a\}$ is a cyclic group of order $2$, 
and where the entries of $P$ are given by $p_{ii}=a$, $p_{ij}=e$ ($i,j\in\mathbb{N}$, $i\neq j$).
Let $M=S^1$ be the monoid obtained from $S$ by adjoining an identity $1$.
Clearly, the rank of $P$ is infinite, and so, by Theorem \ref{thm8}, $M$ is not residually finite.
The maximal subgroups of $M$ are trivial or cyclic of order 2, so finite and hence trivially residually finite.
The $\mathcal{L}$-classes of $M$
are $\{ 1\}$ and $L_j=\{ (i,g,j)\::\: i\in I,\ g\in G\}$ ($j\in J$).
The action of $M$ on $M/\mathcal{L}$ is very simple: $\{1\} (k,g,l)=L_j(k,g,l)=L_l$.
It follows easily that every equivalence relation on $M$ that contains $\mathcal{L}$ is a right congruence,
which implies that the action of $M$ on $M/\mathcal{L}$ is residually finite.
An analogous argument proves that the left action of $M$ on $M/\mathcal{R}$ is residually finite.
\end{proof}

\section{Monoids with $\mathcal{J}$-classes of Finite Type}
\label{sec4}

We will consider the same questions as in the last section, but we will impose a restriction on our monoid $M$ that all its $\mathcal{J}$-classes have finite type, i.e. have finitely many $\mathcal{R}$- and $\mathcal{L}$-classes.
Let us this time start from the regular case, because here we have the following positive result, the second part of which is due to Golubov \cite{Golubov1}.
We note that the construction we use in the proof below also arises in algebraic automata theory where it is known as the \emph{right letter mapping}; see \cite[Chapter~8]{Arbib}.

\begin{theorem}
\label{thm10}
Let $M$ be a regular monoid in which every $\mathcal{J}$-class has finitely many $\mathcal{R}$- and $\mathcal{L}$-classes.
\begin{enumerate}
\item[\textup{(i)}]
The actions of $M$ on its $\mathcal{R}$- and $\mathcal{L}$-classes are residually finite.
\item[\textup{(ii)}]
$M$ is residually finite if and only if all its maximal subgroups are residually finite.
\end{enumerate}
\end{theorem}

\begin{proof}
(i)
Let $a,b\in M$ be any two elements such that $(a,b)\not\in\mathcal{L}$.
We need to show that there exists a finite index right congruence $\rho$ such that
$\mathcal{L}\subseteq \rho$ and $(a,b)\not\in\rho$.
Recall that $M$ acts on $M/\mathcal{L}$ via $L_x \cdot s=L_{xs}$ ($x,s\in M$).

For any $\mathcal{L}$-class $L$ of $M$, let $\mathcal{C}(L)$ be the \emph{cone} of $L$ under the action of $M$:
$$
\mathcal{C}(L)=\{ L\cdot s \::\: s\in M\},
$$
and let $\mathcal{O}(L)$ be the (strong) \emph{orbit} of $L$:
$$
\mathcal{O}(L)=\{ L^\prime \in M/\mathcal{L} \::\: L^\prime\in\mathcal{C}(L)\ \&\ L\in\mathcal{C}(L^\prime)\}.
$$
Note that $\bigcup\mathcal{O}(L)$ is the $\mathcal{D}$-class of $M$ containing $L$.
The assumption that there are finitely many $\mathcal{L}$-classes in every $\mathcal{J}$-class
implies that all $\mathcal{O}(L)$ are finite.
It is clear that the action of $M$ on $M/\mathcal{L}$ induces an action on $\mathcal{C}(L)$.
Moreover, the set $\mathcal{C}(L)\setminus\mathcal{O}(L)$ is an \emph{ideal} in this action:
$$
L^\prime \in \mathcal{C}(L)\setminus\mathcal{O}(L)\ \&\ s\in M \Rightarrow L^\prime\cdot s\in
\mathcal{C}(L)\setminus\mathcal{O}(L).
$$
Therefore, there is an induced action of $M$ on $\mathcal{O}(L)\cup \{0\}$, where $L^\prime\cdot s=0$ whenever
the resulting $\mathcal{L}$-class is in $\mathcal{C}(L)\setminus\mathcal{O}(L)$.
In particular, for every $s\in M$ there exist mappings
$$
\alpha_s\::\: \mathcal{O}(L_a)\cup\{0\}\rightarrow\mathcal{O}(L_a)\cup\{0\},\ 
\beta_s\::\: \mathcal{O}(L_b)\cup\{0\}\rightarrow\mathcal{O}(L_b)\cup\{0\},
$$
where, for $L\in\mathcal{O}(L_a)$ and $L^\prime\in\mathcal{O}(L_b)$, we have
$$
L\alpha_s=\left\{ \begin{array}{ll} L\cdot s & \mbox{if } L\cdot s\in\mathcal{O}(L_a)\\
                                    0 & \mbox{otherwise}
                  \end{array} \right. ,\ 
L^\prime \beta_s=\left\{ \begin{array}{ll} L^\prime \cdot s & \mbox{if } L^\prime\cdot s\in\mathcal{O}(L_b)\\
                                    0 & \mbox{otherwise}
                  \end{array} \right. .
$$
Furthermore, let
$$
A(s)=\{ \alpha_{us}\::\: u\in M\},\ B(s)=\{ \beta_{us}\::\: u\in M\},
$$
and define a relation $\rho$ on $M$ by:
$$
\rho=\{ (x,y)\::\: A(x)=A(y)\ \&\ B(x)=B(y)\}.
$$
It is clear that $\rho$ is an equivalence relation.
Also, as $\mathcal{O}(L_a)$ and $\mathcal{O}(L_b)$ are finite, it follows that the sets $\{ \alpha_s\::\: s\in M\}$ and
$\{\beta_s\::\: s\in M\}$ are finite, and hence so are the sets $\{ A(s)\::\: s\in M\}$ and $\{B(s)\::\: s\in M\}$.
Therefore, the equivalence $\rho$ has finite index.

We now claim that $\rho$ is a right congruence.
To prove this, suppose we have $(x,y)\in\rho$ and $z\in M$.
Consider an arbitrary $\alpha_{uxz}\in A(xz)$.
From $(x,y)\in\rho$ it follows that $\alpha_{ux}=\alpha_{vy}$ for some $v\in M$.
Since the mappings $\alpha_s$ arise from an action of $M$, we have
$$
\alpha_{uxz}=\alpha_{ux}\alpha_z=\alpha_{vy}\alpha_z=\alpha_{vyz}\in A(yz).
$$
This shows that $A(xz)\subseteq A(yz)$. By symmetry we have $A(yz)\subseteq A(xz)$ and so $A(xz)=A(yz)$.
An analogous argument shows that $B(xz)=B(yz)$ and therefore $(xz,yz)\in\rho$, proving that $\rho$ is a right congruence.

Finally, we claim that $(a,b)\not\in\rho$.
Suppose to the contrary.
Since $(a,b)\not\in\mathcal{L}$ without loss of generality we may suppose that $b\not\in Ma$.

We claim that for any $L\in\mathcal{O}(L_b)$ we either have $L\cdot b=L_b$ or else $L\cdot b\not\in \mathcal{O}(L_b)$.
Indeed, finite $\mathcal{J}$-type of $S$ implies the minimum condition on $\mathcal{L}$-classes in any $\mathcal{J}$-class of $S$, and it follows by \cite[Lemma 6.41]{CliffordAndPreston} that $L_b$ is minimal in its $\mathcal{J}$-class $J_b$.
Suppose $L=L_d$, so that $L\cdot b=L_{db}\leq L_b$.
If $L\cdot b\in\mathcal{O}(L_b)$ then $L\cdot b\subseteq J_b$, and so it follows from minimality of $L_b$ that $L\cdot b=L_b$.
 
Furthermore, since $M$ is regular, the $\mathcal{R}$-class of $b$ contains an idempotent $e$, and they satisfy $eb=b$,
which implies that $L_e\cdot b=L_b$. We can conclude that $L_b\in \im(\beta_b)$.
Since we supposed that $B(a)=B(b)$, there must exist $u\in M$ such that $L_b\in \im(\beta_{ua})$.
So for some $L\in\mathcal{O}(L_b)$ we have $L\beta_{ua}=L_b$. 
Thus, for some $c\in L$ we have $cua=b$, which contradicts $b\not\in Ma$, completing the proof of this part.
\vspc

(ii)
This is the main theorem in \cite{Golubov1}.
\end{proof}

Sadly, for non-regular semigroups the situation is not at all so nice.
In order to exhibit examples which show that the analogue of Theorem \ref{thm10} (ii) ($\Leftarrow$) fails
we utilise the following:

\begin{theorem}[Gray, Ru\v{s}kuc \cite{Gray2}]
\label{thmrfdp}
The direct product $S \times T$ of two semigroups is residually finite if and only if $S$ and $T$ are both residually finite. 
\end{theorem}

\begin{example}
\label{ex12a}
Let $S$ be be any semigroup with the property that 
\begin{equation}
\label{eq12b}
xy\not\in\{x,y\}\ (x,y\in S);
\end{equation}
for instance we may take $S=\mathbb{N}$, the additive semigroup of natural numbers.
Let $T$ be any non-residually-finite group $G$; e.g. $T=\mathbb{Q}$, the additive group of rationals, will do.
Let $M$ be the direct product $S\times T$ with an identity adjoined to it.
By Theorem \ref{thmrfdp} we have that $M$ is not residually finite.
On the other hand, property (\ref{eq12b}) implies that $\mathcal{J}=\mathcal{R}=\mathcal{L}=\mathcal{H}=\Delta_M$, and so its $\mathcal{J}$-classes have finite type, and all its Sch\"{u}tzenberger groups are trivial.
\end{example}

Since $\mathcal{R}$ is trivial in the above example, it follows that the action of $M$ on its $\mathcal{R}$-classes is not residually finite. In the following example we exhibit a monoid $M$ in which $\mathcal{J}$-classes have finite type, 
all Sch\"{u}tzenberger groups are residually finite, the actions on $\mathcal{R}$- and $\mathcal{L}$-classes are residually finite, but $M$ is still not residually finite!

\begin{example}
\label{ex13}
We will define $M$ by means of a presentation in stages.
The generators are
\begin{equation}
\label{eq1}
a_i\ (i\in \mathbb{Z}),\ b_j\ (j\in\mathbb{Z}),\ c_k\ (k\in\mathbb{N}),\ d.
\end{equation}

The generators $a_i$ generate a free semigroup $A$ and the generators $c_k$ generate a free semigroup $C$.
The generators $b_j$ serve as reference points for a sequence $B_j$ of $\mathcal{H}$-classes, which also turn out to be
 $\mathcal{R}$- and $\mathcal{L}$-classes. The Sch\"{u}tzenberger group of $B_j$ will be isomorphic to 
the direct product $C_2\times C_2\times\ldots$ of cyclic groups of order $2$.
On the right this will be generated by all the $c_k$, while on the left it will be generated by the $a_i$ ($i>j$);
the generator $a_i$ on the left corresponds to the generator $c_{i-j}$ on the right; the generators
$a_i$ ($i\leq j$) act trivially:
\begin{eqnarray}
\label{eq4a}
&&a_i b_j =b_j\ (i,j\in\mathbb{Z}, i\leq j);\\
\label{eq4b}
&&a_i^2 b_j= b_j \ (i,j\in\mathbb{Z}, i>j);\\
\label{eq4c}
&&a_ia_lb_j=a_la_ib_j\ (i,l,j\in\mathbb{Z});\\
\label{eq4d}
&&b_jc_k^2=b_j\ (j\in\mathbb{Z},\ k\in\mathbb{N});\\
\label{eq4e}
&&b_jc_kc_m=b_jc_mc_k\ (j\in\mathbb{Z},\ k,m\in\mathbb{N});\\
\label{eq4f}
&&a_i b_j =  b_j c_{i-j} \ (i,j\in\mathbb{Z},\ i>j).
\end{eqnarray}
(Note that the relations (\ref{eq4d}) and (\ref{eq4e}) are in fact redundant.)

The generator $d$ generates a free monogenic semigroup.
Multiplication by $d$ on the right will take $B_j$ to $B_{j+1}$:
\begin{equation}
\label{eq5}
b_jd=b_{j+1} \ (j\in\mathbb{Z}).
\end{equation}
Finally, we have the following zero products:
\begin{eqnarray}
a_ic_k = c_ka_i=b_ja_i  = c_kb_j & = & b_jb_l=a_id \notag \\
 & = &  da_i=dc_k =  db_j=0\ (i,j,l\in\mathbb{Z},\ k\in\mathbb{N}).
\end{eqnarray}

A set of normal forms for $M$ is
\begin{equation}
\label{eq6}
M=\{1\}\cup A\cup C \cup (\bigcup_{j\in\mathbb{Z}} B_j) \cup D \cup \{0\}
\end{equation}
where
\begin{eqnarray*}
&&B_j=\{ b_j c_{k_1} c_{k_2}\ldots c_{k_t}\::\:
t\geq 0,\ k_1<k_2<\ldots<k_t\}\ (j\in\mathbb{Z}),\\
&& D=\{ wd^t\::\: w\in C^1,\ t\in\mathbb{N}\}.
\end{eqnarray*}

In this monoid all Green's equivalences coincide; 
the only non-singleton classes are
$B_j$ ($j\in\mathbb{Z}$).
Trivially, $M$ has finite $\mathcal{J}$-type. 
The Sch\"{u}tzenberger groups of trivial $\mathcal{H}$-classes are trivial. The Sch\"{u}tzenberger group of $B_j$ is $C_2\times C_2\times \ldots$. In particular, all the Sch\"{u}tzenberger groups are residually finite.

The actions of $M$ on its $\mathcal{L}$- and $\mathcal{R}$-classes are as follows:
\begin{gather*}
\begin{array}{r|cccc}
S/\mathcal{L} & a_i & c_k & b_j & d\\ \hline
1 & a_i & c_k & B_j & d\\
w\in A & wa_i & 0 & B_j & 0\\
w\in C & 0 & wc_k & 0 & wd\\
wd^t\in D & 0 & wc_kd^t & 0 &  wd^{t+1}\\
B_j & 0 & B_j & 0 & B_{j+1}\\
0 & 0 & 0 & 0 & 0
\end{array}
\\
\begin{array}{ r|cccccc}
S/\mathcal{R} & 1 & w\in A & w\in C & wd^t\in D & B_j & 0\\ \hline
a_i & a_i & a_iw & 0 & 0 & B_j &0 \\
c_k & c_k & 0 & c_kw & c_k wd^t & 0 & 0 \\
b_j & B_j & 0 & B_j & B_{j+t} & 0 & 0\\
d & d & 0 & wd & wd^{t+1} & 0 & 0
\end{array}
\end{gather*}
and are sketched in Figures \ref{fig1}
 and \ref{fig2} respectively.

\setcounter{figure}{0}
\renewcommand{\figurename}{Figure}

\begin{figure}
\begin{center}
\includegraphics{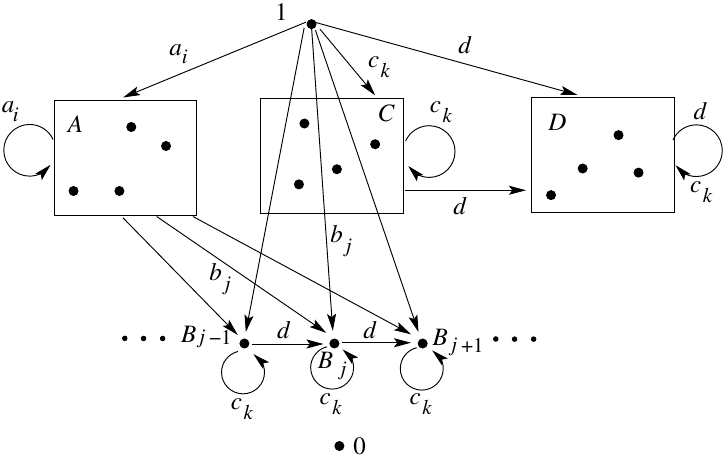}
\end{center}
\caption{The action of $M$ on its $\mathcal{L}$-classes. The missing arrows point to $0$.}
\label{fig1}
\end{figure}

\begin{figure}
\begin{center}
\includegraphics{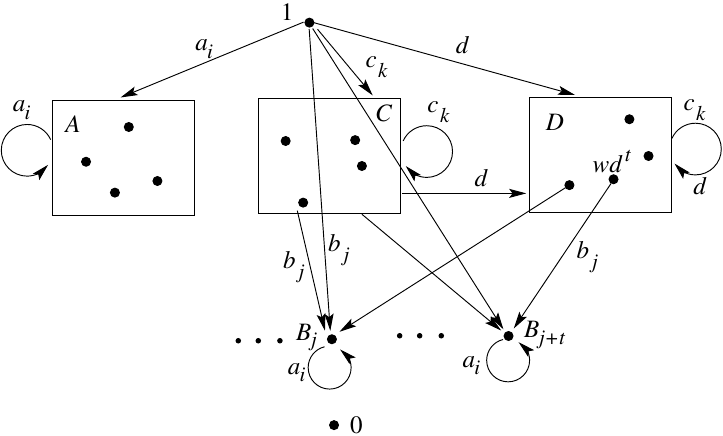}
\end{center}
\caption{The left action of $M$ on its $\mathcal{R}$-classes. The missing arrows point to $0$.}
\label{fig2}
\end{figure}

These two actions are residually finite.
Indeed, if we want to separate two distinct elements $u,v\in M$ such that at least one of them belongs to the set $\{1\}\cup A\cup C\cup D$, then we can let $N=\max(|u|,|v|)$ and
$$W=\{w\in \{1\}\cup A\cup C\cup D\::\: |w|\leq N\},$$ 
and then prove that the relation $\Delta_W\cup\Phi_{S\setminus W}$ is a congruence for both actions, and that it separates $u$ and $v$.
In order to separate $B_j$ and $B_l$ ($j\neq l$) we can choose $M\in\mathbb{N}$ such that
$j\not\equiv l\pmod{M}$, and use the relation
\begin{eqnarray*}
&&\Phi_A\cup\Phi_C \\ &\cup& \{ (w_1d^{t_1},w_2d^{t_2})\::\: w_1,w_2\in C^1,\ t_1,t_2\in\mathbb{N},\ 
t_1\equiv t_2 \pmod{M}\}\\
&\cup&
\{ (B_{m_1},B_{m_2})\::\: m_1\equiv m_2\pmod{M}\}\cup\{(0,0)\}.
\end{eqnarray*}
Finally, the same relation, with $M$ chosen arbitrarily, will separate any element of any $B_j$ from $0$.

Finally, we prove that $M$ is not residually finite. Suppose $f\::\: M\rightarrow N$ is a homomorphism, where $N$ is finite.
We shall prove that
\begin{equation}
\label{eq12}
f(b_0c_1)=f(b_0).
\end{equation}
From the infinite sequence $c_{1},c_{2},c_{3},\ldots$ pick a pair of entries
$c_{p}$, $c_{q}$ such that
\begin{equation}
\label{eq13}
f(c_{p})=f(c_{q}) \mbox{ and } 1\leq p<q.
\end{equation}
Such a pair must exist because $N$ is finite.
Multiply both sides of (\ref{eq13}) by $f(b_{1-q})$ on the left:
\begin{equation}
\label{eq14}
f(b_{1-q}c_{p})=f(b_{1-q}c_{q}).
\end{equation}
Using relations (\ref{eq4f}) transform (\ref{eq14}) into
\begin{equation}
\label{eq15}
f(a_{1+p-q}b_{1-q})=f(a_{1}b_{1-q}).
\end{equation}
Now multiply both sides by $f(d^{q-1})$ on the right:
\begin{equation}
\label{eq16}
f(a_{1+p-q}b_{1-q}d^{q-1})=f(a_{1}b_{1-q}d^{q-1}),
\end{equation}
and use relations (\ref{eq5}):
\begin{equation}
\label{eq17}
f(a_{1+p-q}b_{0})=f(a_{1}b_{0}).
\end{equation}
By (\ref{eq4f}) we have $a_{1}b_0=b_0 c_1$. Since $1+p-q\leq 0$, by
(\ref{eq4a}) we have $a_{1+p-q}b_{0}=b_0$. Thus (\ref{eq17}) implies
$f(b_0c_1)=f(b_0)$, as required. This completes Example \ref{ex13}.
\end{example}

As for the other direction, we will postpone its consideration until Section \ref{sec4b}.

\section{Intermezzo: Residual Finiteness, Ideals and Rees Quotients}
\label{sec4a}

If $I$ is an ideal of a monoid $M$ then the relation $\Delta_M\cup\Phi_I$ is easily seen to be a congruence on $M$, and the corresponding quotient is called the \emph{Rees quotient} of $M$ by $I$ and denoted by $M/I$.
It can be identified with the set $(M\setminus I)\cup\{0\}$, where the multiplication is as in $M$, except that all the elements of $I$ are identified with the zero element $0$.
In this section we will prove a technical result and its corollary which describe a couple of situations where residual finiteness of a monoid is completely determined by residual finiteness of an ideal and the corresponding Rees quotient.
A more detailed analysis of Rees quotients and residual finiteness can be found in the article \cite{GolubovSapir} by Golubov and Sapir.
One result they prove is that if $I$ is an ideal of a semigroup $S$ such that $S\setminus I$ is finite
then $S$ is residually finite if and only if $I$ is residually finite.
A more general version of this that we will use below is:

\begin{theorem}[Ru\v{s}kuc, Thomas \cite{Ruskuc&Thomas}]
\label{thm4a1a}
Let $S$ be a semigroup, and let $T$ be a subsemigroup of $S$ with $S\setminus T$ finite.
Then $S$ is residually finite if and only if $T$ is residually finite.
\end{theorem}

We begin with a definition:

\begin{definition}
Let $M$ be a monoid, and let $I$ be an ideal of $M$. We say that $I$ is \emph{rf-compatible} with $M$ if for any two distinct $s,t\in I$ there exists a finite index congruence $\rho$ on $I$ such that $s/\rho\neq t/\rho$ and the relation $\rho\cup\Delta_M$ is a congruence on $M$.
\end{definition}

\begin{remark}
We note that rf-compatibility of $I$ implies its residual finiteness.
\end{remark}

\begin{proposition}
\label{prop4a1}
Let $M$ be a monoid with zero and let $I$ be an rf-compatible ideal of $M$ such that 
$N=M\setminus I \cup\{0\}$ is a submonoid of $M$.
Then $M$ is residually finite if and only if $N$ is residually finite.
\end{proposition}

\begin{proof}
($\Rightarrow$)
This is obvious, since $N$ is a submonoid of $M$.

($\Leftarrow$)
Let $s,t\in M$ be two any two distinct elements. 
Consider the set $A=\{s,t,0\}\cap I$.
Since $I$ is rf-compatible, there exists a finite index congruence $\rho$ on $I$ which separates any distinct elements
in $A$ (when $|I|=1$ we can freely choose $\rho$) and $\sigma=\rho\cup\Delta_M$ is a congruence on $M$.
The quotient $Q=M/\sigma$ contains (a copy of) $N$, and $Q\setminus N$ is finite.
By Theorem \ref{thm4a1a}, $Q$ is residually finite since $N$ is residually finite.
Since $s/\sigma\neq t/\sigma$, we can separate $s/\sigma$ and $t/\sigma$ by a finite quotient of $Q$,
which is also a finite quotient of $M$. \end{proof}

\begin{corollary}
\label{corol4a5}
Let $M$ be a monoid, and let $I$ be an rf-compatible ideal of $M$ such that 
$N=M\setminus I$ is a submonoid of $M$.
Then $M$ is residually finite if and only if $N$ is residually finite.
\end{corollary}

\begin{proof}
Adjoin a zero element to $M$ (and to $I$ and $N$), and use Proposition \ref{prop4a1}.
\end{proof}

Golubov \cite{Golubov1} (immediately following Corollary 3) exhibits an example of a residually finite inverse semigroup $S$ which has a non-residually finite Rees quotient $S/I$.
The following example shows that
the following statement is also not true: if $I$ and $M/I$ are residually finite then $M$ is residually finite.

\begin{example}
Let $G$ be a residually finite group acting primitively on an infinite set $X$, via $(x,g)\mapsto x\cdot g$.
(For a specific example, one could take $G$ to be a free group of countably infinite rank, and let it act on a countable set $X$ via the epimorphism onto the finitary symmetric group on $X$.)
Define a monoid
$$
M=\{1\}\cup G\cup X\cup\{0\},
$$
where the non-obvious part of multiplication is defined by:
$$
xg=x\cdot g,\ gx=xy=0\ (g\in G,\ x,y\in X).
$$
Clearly, the set $I=X\cup \{0\}$ is an ideal with zero multiplication; 
in particular, it is residually finite.
The quotient $M/ I$ is isomorphic to $\{1\}\cup G\cup\{0\}$, and is clearly residually finite, since $G$ is.
However, $M$ itself is not residually finite.
Indeed, if $\rho$ is any congruence on $M$, the restriction $\rho\restriction_X$
is a $G$-invariant partition of $X$, which, by primitivity, has to be $\Delta_X$ or $\Phi_X$.
It follows that finite index congruences cannot separate pairs of elements from $X$.
\end{example}

\section{Monoids with $\mathcal{J}$-classes of Finite Type (continued)}
\label{sec4b}

We complete the consideration of monoids with $\mathcal{J}$-classes of finite type by exhibiting such a monoid which is residually finite, but in which the actions on $\mathcal{R}$- and $\mathcal{L}$-classes are not residually finite.

Before we  write down the presentation for our monoid we need to introduce an auxiliary function
$\tau : \mathbb{Z}\setminus\{0\}\rightarrow\mathbb{Z}$:
\begin{equation}
\label{eq4b1}
\tau(2^k(2r+1))= \frac{2}{3} (2^{2\lceil k/2\rceil}-1)\ (k,r\in\mathbb{Z},\ k\geq 0),
\end{equation}
where $\lceil x\rceil$ denotes the smallest integer not smaller than $x$.
We record two properties of $\tau$ that will be of use in what follows:

\begin{lemma}
\label{lemma4b6}
Let $a,b\in\mathbb{Z}\setminus\{0\}$ and $m\in\mathbb{N}$.
If $a\equiv b \pmod{2^m}$ then $\tau(a)\equiv \tau(b) \pmod{2^{m+1}}$.
\end{lemma}

\begin{proof}
Suppose $a=2^k(2r+1)$, $b=2^l(2u+1)$. From $a\equiv b\pmod{2^m}$ it follows that $k=l$ or else $k,l\geq m$.
If $k=l$ then $\tau(a)=\tau(b)$. Otherwise, we have $2\lceil k/2\rceil , 2\lceil l/2 \rceil \geq m$
and so
$$
\tau(a)-\tau(b)= \frac{2}{3}(2^{2\lceil k/2\rceil}-2^{2\lceil l/2\rceil})
$$
is divisible by $2^{m+1}$, as required.
\end{proof}

\begin{lemma}
\label{lemma4b6a}
There is no $x\in\mathbb{Z}$ satisfying $x\equiv \tau(2^m) \pmod{2^{m+1}}$ for all $m=1,2,\ldots$.
\end{lemma}

\begin{proof}
This follows immediately from $0<\tau(2^m) < 2^{m+1}$, and the fact that both
$\tau(2^m)$ and
$$
2^{m+1}-\tau(2^m)= \frac{2}{3}( 2^{2\lceil (m-1)/2 \rceil +1} +1)
$$
are unbounded as $m$ increases.
\end{proof}

We are now ready to begin defining our monoid $M$. It is a commutative monoid with zero, with generators
\begin{equation}
\label{eq4b2}
\{ a,a^{-1}\} \cup \{ b_i,c_i\::\: i\in\mathbb{Z}\}\cup \{d,e\}.
\end{equation}
It has an infinite cyclic group of units generated by $a$:
\begin{equation}
\label{eq4b3}
aa^{-1}=a^{-1}a=1.
\end{equation}
The multiplication of $b_i$ and $c_j$ ($i,j\in\mathbb{Z}$) is governed by
\begin{equation}
\label{eq4b4}
b_ic_j = \left\{ \begin{array}{ll} 
                   d & \mbox{if } i=j\\
                   a^{\tau(j-i)}e & \mbox{if } i\neq j
                 \end{array}\right.
\end{equation}
and all other products of two generators are zero:
\begin{equation}
\label{eq4b5}
b_ib_j=b_ic_j=b_id=b_ie=c_jc_k=c_jd=c_je=dd=de=ee=0\ (i,j,k\in\mathbb{Z}).
\end{equation}
A standard argument shows that the following is a set of distinct normal forms for $M$:
\begin{equation}
\label{eq4b8}
M=A\cup B\cup C\cup D\cup E\cup\{0\},
\end{equation}
where
\begin{eqnarray}
\label{eq4b9}
&& A=\{a^{\pm p} \::\: p\in\mathbb{Z}\},\\
\label{eq4b10}
&&
B=\displaystyle\bigcup_{i\in\mathbb{Z}} B_i,\ B_i=Ab_i\ (i\in\mathbb{Z}),\\
\label{eq4b11}
&&
C=\displaystyle\bigcup_{i\in\mathbb{Z}} C_i,\ C_i=Ac_i\ (i\in\mathbb{Z}),\\
\label{eq4b12}
&&
D=Ad,\ E=Ae.
\end{eqnarray}
Since $M$ is commutative it follows that all Green's relations coincide, and, in particular, $\mathcal{J}$-classes are of finite type.
Also, from the above normal forms, it is immediate to see that the equivalence classes of any Green's relation are precisely the sets $A$, $B_i$, $C_i$ ($i\in I$), $D$, $E$ and $\{0\}$.
The (left- and right-) action of $M$ on these classes and is shown in Table~\ref{tab4b13a} and illustrated in Figure \ref{fig4b13}.

\setcounter{figure}{1}
\renewcommand{\figurename}{Table}

\begin{figure}
$$
\begin{array}{c|ccccc}
  & a^{\pm 1} & b_i & c_i & d & e \\ \hline
A & A & B_i & C_j & D & E \\
B_i & B_i & 0 & D &  0& 0 \\
B_j & B_j &  0 & E  & 0&0\\ 
C_i & C_i & D & 0&0&0\\
C_j&C_j & E & 0&0&0 \\
D&D&0&0&0&0\\
E&E&0&0&0&0\\
0&0&0&0&0&0
\end{array}
$$
\caption{The action of $M$ on its Green's classes (where $i\neq j$).}
\label{tab4b13a}
\end{figure}

\setcounter{figure}{2}
\renewcommand{\figurename}{Figure}

\begin{figure}
\begin{center}
\includegraphics{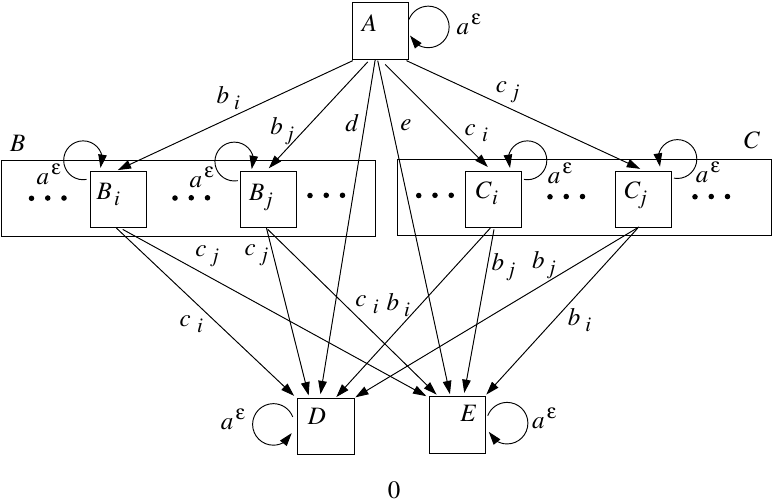}
\caption{The action of $M$ on its Green's classes. The arrows not shown all point to $0$.}
\label{fig4b13}
\end{center}
\end{figure}

\begin{lemma}
The action of the monoid $M$ on its $\mathcal{L}$-classes is not residually finite.
\end{lemma}

\begin{proof}
We are going to show that for any finite index right congruence $\rho$ on $M$ that contains $\mathcal{L}$, the elements
$d$ and $e$ are in the same class (while, clearly, they are not $\mathcal{L}$-equivalent in $M$).
Since $\rho$ has finite index, there exist distinct $i,j\in\mathbb{Z}$ such that $(b_i,b_j)\in\rho$.
But then $\rho\ni (b_ic_j,b_jc_j) = (a^{\tau(j-i)}e,d)$. Since $(a^{\tau(j-i)}e,e)\in\mathcal{L}\subseteq \rho$, it follows that $(d,e)\in\rho$, as claimed.
\end{proof}

We now start on the harder task of proving that $M$ is residually finite.
We will achieve this by twice decomposing it in line with Proposition \ref{prop4a1} and Corollary \ref{corol4a5}.

\begin{lemma}
\label{lemma4b14}
The set $I=C\cup D\cup E\cup\{0\}$ is an rf-compatible ideal of $M$, and the set $N=M\setminus I \cup \{0\}=A\cup B\cup \{0\}$ is a submonoid of $M$.
\end{lemma}

\begin{proof}
It is straightforward to check that $I$ is an ideal, and that $N$ is a submonoid, and it remains 
to prove that $I$ is rf-compatible.
To this end let $s,t\in I$ be any two distinct elements.
We need to find a congruence of finite index $\rho$ on $I$ such that $s/\rho\neq t/\rho$ and
$\rho\cup\Delta_M$ is a congruence on $M$.
We note that $I$ is a zero semigroup, and so any equivalence relation on it is a congruence.
We split our considerations into a number of cases, depending from which part of $I$ the elements $s$ and $t$ come.

\textit{Case 1: $s\in C$,} say $s=a^pc_i$. If $t\in C_i$, then write $t=a^qc_i$, where $p\neq q$,
and 
let $m\in \mathbb{N}$ be such that $m>|p-q|$, so that $p\not\equiv q\pmod{m}$.
Otherwise pick an arbitrary $m\in\mathbb{N}$.
Define $\rho$ by:
$$
\rho=\{ (a^rc_i, a^uc_i)\::\: r\equiv u\pmod{m}\}\cup \Phi_{I\setminus C_i}.
$$
It is clear that $s/\rho\neq t/\rho$. 
To show that $\rho\cup\Delta_M$ is a congruence, it is sufficient to show that for any $(x,y)\in\rho$
and any generator $g$ of $M$ we have $(gx,gy)\in\rho$.
If $x=a^rc_i$, $y=a^u c_i$ (with $r\equiv u\pmod{m}$) and $g=a^\epsilon$ ($\epsilon=\pm 1$) 
then we have
$$
(g x,g y)=(a^{r+\epsilon}c_i, a^{u+\epsilon}c_i)\in \rho,
$$
while for any other possibility for $x$, $y$ and $g$ we have $gx,gy\in I\setminus C_i$, implying
$(gx,gy)\in\Phi_{I\setminus C_i}\subseteq \rho$.

\textit{Case 2: $s,t\in D$,} say $s=a^pd$, $t=a^qd$, $p\neq q$.
Let $m\in \mathbb{N}$ be such that $2^{m+1}>|p-q|$, so that $p\not\equiv q\pmod{2^{m+1}}$.
Define $\rho$ by:
\begin{eqnarray}
\label{eq4b21}
\rho &=&
\{ (a^r c_i, a^u c_j)\::\: r\equiv u\pmod{2^{m+1}},\ i\equiv j\pmod{2^m} \} \\
\label{eq4b22}
&&
\cup \{ (a^r d, a^u d)\::\: r\equiv u\pmod{2^{m+1}} \} \\
\label{eq4b23}
&&
\cup \{ (a^r e, a^u e)\::\: r\equiv u\pmod{2^{m+1}} \} \\
\label{eq4b24}
&&
\cup \{ (a^r d, a^u e), (a^u e,a^r d) \::\: r+\tau(2^m) \equiv u\pmod{2^{m+1}} \} \\
\label{eq4b25}
&&
\cup \{ (0,0)\}.
\end{eqnarray}
To see that $\rho$ is an equivalence relation, observe that the restriction of $\rho$ on each of
$C$, $D$, $E$ and $\{0\}$ is clearly an equivalence relation, and that the only additional relationships serve to identify the equivalence class of $a^rd$ in $D$ with the equivalence class of  $a^{r+\tau(2^m)}e$ in $E$.
The choice of $m$ implies that $s/\rho\neq t/\rho$, and it only remains to be proved that $\rho\cup\Delta_M$ is a congruence.
Multiplying a pair from any of the sets (\ref{eq4b21})--(\ref{eq4b25}) by $a^\epsilon$ ($\epsilon=\pm 1$)
yields a pair from the same set.
Multiplying such a pair by $c_k$ ($k\in\mathbb{Z}$), $d$ or $e$,
or multiplying a pair from (\ref{eq4b22})--(\ref{eq4b25}) by $b_k$ ($k\in\mathbb{Z}$) always yields $(0,0)\in\rho$.
There remains to analyse the effect of multiplying a typical pair from (\ref{eq4b21}) by $b_k$:
\begin{numcases}{(b_ka^rc_i,b_ka^uc_j)=}
\label{eq4b31}
(a^{r+\tau(i-k)}e,a^{u+\tau(j-k)}e) & if $i\neq k\neq j$\\
\label{eq4b32}
(a^rd,a^{u+\tau(j-k)}e) & if $i=k\neq j$\\
\label{eq4b33}
(a^{r+\tau(i-k)}e,a^ud) & if $i\neq k= j$\\
\label{eq4b34}
(a^rd ,a^ud) & if $i= k= j$.
\end{numcases}
Note that $i-k\equiv j-k \pmod{2^m}$ and so, by Lemma \ref{lemma4b6}, we have
$\tau(i-k)\equiv \tau(j-k)\pmod{2^{m+1}}$. Since, in addition, $r\equiv u\pmod{2^{m+1}}$, it follows
that pair (\ref{eq4b31}) belongs to set (\ref{eq4b23}).
Now let us consider pair (\ref{eq4b32}). From $i\equiv j\pmod{2^m}$ and $k=i$ we have
$j-k\equiv 0\equiv 2^m \pmod{2^m}$, and so, by Lemma \ref{lemma4b6}, we conclude that
$\tau(j-k)\equiv \tau(2^m)\pmod{2^{m+1}}$. Combining this with $r\equiv u\pmod{2^{m+1}}$, we deduce
that pair (\ref{eq4b32}) belongs to set (\ref{eq4b24}).
The pair (\ref{eq4b33}) is symmetric to (\ref{eq4b32}) and also belongs to set (\ref{eq4b24}).
Finally, pair (\ref{eq4b34}) clearly belongs to (\ref{eq4b22}), completing Case 2.

\textit{Case 3: $s\in D$, $t\in E$,} say $s=a^pd$, $t=a^qe$.
Pick $m\in\mathbb{N}$ such that
\begin{equation}
\label{eq4b41}
p+\tau(2^m)\not\equiv q\pmod{2^{m+1}};
\end{equation}
the existence of such $m$ is guaranteed by Lemma \ref{lemma4b6a}.
Define a relation $\rho$ by (\ref{eq4b21})--(\ref{eq4b25}) above.
The above proof that $\rho\cup\Delta_M$ is a congruence on $M$ remains valid,
and (\ref{eq4b41}) implies that $s/\rho\neq t/\rho$.

\textit{Case 4: $s\in D\cup E$, $t=0$.}
Note that the congruence $\rho$ constructed in Case 2 (with an arbitrary choice of $m$)
has the property $0/\rho=\{0\}$, and so separates $s$ and $t$.

\textit{Case 5: $s, t\in E$,} say $s=a^pe$, $t=a^qe$, $p\neq q$.
The congruence $\rho$ and the proof carry over from Case 3 verbatim,
completing the proof of this case, and of the whole lemma.
\end{proof}

Let us now examine the submonoid $N=M\setminus I \cup\{0\}=A\cup B\cup\{0\}$.

\begin{lemma}
\label{lemma4b15}
In the monoid $N$, the set $B\cup \{0\}$ is an rf-compatible ideal, and the set $A=N\setminus (B\cup\{0\})$
is a submonoid.
\end{lemma}

\begin{proof}
This proof is very similar to, and much easier than, that of Lemma \ref{lemma4b14}.
Again, the only non-obvious assertion is rf-compatibility.
Without loss of generality assume that $s\in B$, say $s=a^p b_i$.
As for $t$, either $t=a^q b_j$ (with $p\neq q$ or $i\neq j$), or else $t=0$, in which case we let $q$ be arbitrary.
Let $m\in \mathbb{N}$ be such that $m>|p-q|$, and let $\rho$ be the relation
$$
\rho = \{ (a^r,a^u),(a^rb_i,a^ub_i)\::\: r\equiv u\pmod{m}\} \cup \Phi_{B\setminus B_i} \cup \{(0,0)\}.
$$
A routine verification shows that $s/\rho\neq t/\rho$ and that $\rho\cup\Delta_N$ is a congruence. \end{proof}

We can now see that the monoid $M$ is residually finite: Firstly, from Lemma \ref{lemma4b15}, Corollary \ref{corol4a5} and the fact that $A$ is the infinite cyclic group, which is residually finite, it follows that $N$ is residually finite. Then from this, Lemma \ref{lemma4b14} and Proposition \ref{prop4a1} it follows that $M$ itself is residually finite. 

To summarise:

\begin{proposition}
\label{prop4b16}
The commutative monoid $M$ with zero defined by generators (\ref{eq4b2}) and defining relations (\ref{eq4b3}), (\ref{eq4b4}), (\ref{eq4b5}) is residually finite, but its action on the $\mathcal{R}$-classes is not residually finite.
\end{proposition}

As a curiosity we mention that the monoid $M$ is very close to being a submonoid of the direct product
of the infinite cyclic group $C_\infty$ and the quotient $M/\mathcal{H}$,
which is residually finite by Theorem \ref{thmrfdp}.
What prevents $M$ from being such a submonoid is part (\ref{eq4b4}) of its multiplication.
Thus, if we chose $\tau$ to be constantly $0$ (or indeed any constant mapping) our example would not work any longer,
since residual finiteness would fail.

\section{Monoids with Finitely Many Left and Right Ideals}
\label{sec5}

We conclude our account by imposing a very strong finiteness condition of there being only finitely many left- and right ideals in our monoid.
Of course, this makes the actions on $\mathcal{L}$- and $\mathcal{R}$-classes finite, and hence residually finite too.
As an immediate consequence of Golubov's Theorem \cite{Golubov1} (see Theorem \ref{thm10} above) we obtain:

\begin{corollary}
\label{corol15}
A regular monoid with finitely many left- and right ideals is residually finite if and only if all its maximal subgroups are residually finite.
\end{corollary}

After the results of the previous sections, where the differences between the regular and non-regular cases were starkly exposed, this time we have a full analogue of the above result:

\begin{theorem}
\label{thm16}
Let $M$ be a monoid with finitely many left and right ideals. Then $M$ is residually finite if and only if all its 
Sch\"{u}tzenberger groups are residually finite.
\end{theorem}

\begin{proof}
The direct part follows from the more general Theorem \ref{thm1}.
For the converse part, 
we need some general theory.
Given an equivalence relation $\pi$ on a monoid $M$ we define 
\[
\Sigma_r( \pi ) = \{ (x,y) \in M \times M : (xm,ym) \in \pi \mbox{ for all } m \in M \}.
\]
Clearly $\Sigma_r( \pi )$ is a right congruence, and it follows from  \cite[Lemma~10.3]{CliffordAndPreston} that $\Sigma_r( \pi )$ is the largest right congruence of $M$ contained in $\pi$. 
The following alternative description of $\Sigma_r( \pi )$ from \cite[Proposition~2.2]{Ruskuc&Thomas} will be useful here:
for $s \in M$ and $X \subseteq M$  define
\[
Q_M(s,X) = \{ x \in M : sx \in X  \};
\]
if $C_i\ (i\in I)$ are the classes of $\pi$, then
\begin{equation}
\label{QS}
(x,y) \in \Sigma_r(\pi) \Leftrightarrow Q_M(x,C_i) = Q_M(y,C_i) \quad \mbox{for all} \ i \in I.
\end{equation}
Analogously we can define $\Sigma_l(\pi)$, the largest left congruence of $S$ contained in $\pi$, and $\Sigma(\pi)$, the largest two-sided congruence of $S$ contained in $\pi$. 

Now let $x,y \in M$ be distinct. 
We need to find a finite index congruence on $M$ separating $x$ and $y$. 
By \cite[Theorem~2.4]{Ruskuc&Thomas}, if $\rho$ is a finite index right congruence then
the two-sided congruence $\Sigma(\rho)\subseteq \rho$ also has finite index,
and so it is in fact sufficient to find a right (or left) congruence on $M$ with finite index that separates $x$ and $y$. 

If $(x,y) \not\in \gl$ then, since $M$ has finitely many $\gl$-classes, it follows that $\gl$ itself is a finite index right congruence separating $x$ and $y$. A dual argument deals with the case $(x,y) \not\in \gr$. 
So from now on we deal with the case where $x\gr y$ and $x\gl y$, i.e. where $x$ and $y$ belong to the same 
$\gh$-class $H$.
Let us fix an element $h \in H$, and let $s_x, s_y \in \mathrm{Stab}(H)$ satisfy $hs_x = x$, and $hs_y = y$ (such elements exist by Proposition~\ref{schutzstabproperties} (iii)). Since the Sch\"{u}tzenberger group $\Gamma(H)=\mathrm{Stab}(H)/\sigma$ is residually finite, there is a normal subgroup $N$ of $\Gamma(H)$ of finite index such that $s_x/\sigma$ and $s_y/\sigma$ belong to different cosets of $N$ in $\Gamma(H)$. 

Let $N_i$ ($i =0,\ldots,m$) be the cosets of $N$ in $\Gamma(H)$ where $N_0 = N$, and let
\[
\overline{N_i} = \{ s \in \mathrm{Stab}(H): s/\sigma \in N_i  \}\ (i=0,\ldots,m).
\]     
Now partition $H$ as $H = \bigcup_{0 \leq i \leq m}{C_i}$ where $C_i = h \overline{N_i}$. Observe that the $C_i$ blocks are preserved by right multiplication from $\mathrm{Stab}(H)$, for
\[
us \in C_j \Leftrightarrow s/\sigma \in {N_i}^{-1}N_j \Leftrightarrow vs \in C_j
\] 
for all $u,v \in C_i$, $s \in M$. Let $\pi$ be the equivalence relation on $M$ with equivalence classes given by the partition $M = C_0 \cup C_1 \cup \ldots C_m \cup C_{m+1}$, where 
$C_{m+1}=M \setminus H $. We claim that the right congruence $\Sigma_r(\pi)$ has finite index. Once established, this will complete the proof of the theorem, since $x$ and $y$ belong to different equivalence classes of 
$\pi$ and hence belong to different classes of $\Sigma_r(\pi)$. 
In order to prove that $\Sigma_r(\pi)$ has finite index it suffices to prove that the collection of sets $Q_M(s,C_k)$, $s \in M$, $0\leq i\leq m+1$ (see \eqref{QS}) is finite. 
Note that $Q_M(s,C_{m+1})=M\setminus (\bigcup_{k=0}^m Q_M(s,C_k))$. Thus it is sufficient to prove the following

\bigskip

\noindent \textbf{Claim.} For every $k=0,\ldots,m$, there are only finitely many different sets $Q_M(s,C_k)$ for $s \in M$.

\bigskip

To prove the claim, let $s \in M$ and fix a finite cross-section of the $\gl$-classes of $M$. 
Let $l$ be the representative of the $\gl$-class of $s$. 
Pick an element $u \in Q_M(s,C_k)$. 
Denote by $L_0$ the $\gl$-class of $M$ that contains $H$. 
Since $\gl$ is a right congruence, we have $su \gl lu$. 
Since $su \in C_k\subseteq H\subseteq L_0$ it follows that $lu \in  L_0$;
let $H'\subseteq L_0$ be the $\gh$-class of $lu$. 
Since $H,H^\prime\subseteq L_0$ it follows from Proposition~\ref{schutzstabproperties} (iv) that $\mathrm{Stab}(H') = \mathrm{Stab}(H)$,
and we can partition $H'$ into blocks $C_i' = h'\overline{N_i}$ where $h'$ is a fixed element of $H'$. 
These blocks have the property that for all $j$, and all $w \in M$, we have $C_j' w = C_j'$ if and only if $w \in \overline{N}$. 
Suppose that $lu \in C_j'$. We claim that
\begin{equation}\label{equation}
Q_M(s,C_k)=\{ u' \in M : su' \in C_k  \}  =  \{ u' \in M: lu' \in C_j'  \}.
\end{equation}    
Indeed, if $su' \in C_k$ then since $su \in C_k$ we can write $su' = suu_2$ for some $u_2 \in M$. Now $su \in C_k$ and $(su)u_2 \in C_k$ which implies that $C_k u_2 = C_k$. It follows that $u_2 \in \overline{N}$ and therefore that $C_j' u_2 = C_j'$. Since $l \gl s$ we can write $l = u_3 s$ where $u_3 \in M$. Now we have
\[
l u' = u_3 s u' = u_3 s u u_2 = (l u) u_2 \in C_j' u_2 = C_j' 
\] 
proving the direct inclusion of \eqref{equation}. The converse inclusion may be proved using a symmetric argument.

Note that the collection of sets $\{ u' \in M: lu' \in C_j'\}$, appearing as the right hand side of \eqref{equation}, is finite. Indeed, it depends only on:
(i) $l$, one of the finitely many $\gl$-class representatives;
(ii) $H^\prime$, one of the finitely many $\gh$-classes in $L_0$; and
(iii) $C_j^\prime$, one of the finitely many blocks of $H^\prime$.
This completes the proof of Claim, and also of the theorem.
\end{proof}

\begin{remark}
Theorem \ref{thm16} in some sense complements
\cite{Ruskuc2} where the analogous statements are shown to hold for the finiteness conditions of being finitely generated, and being finitely presented. 
\end{remark}

\begin{remark}
We observe that Theorem \ref{thm16} and its proof remain valid if finiteness condition on $M$ is weakened to the following:
$M$ has finite $\mathcal{J}$-type, and for every $\mathcal{J}$-class $J$ there are only finitely many $\mathcal{J}$-classes above it.
(Compare this with \cite[Corollary 4.3]{Ruskuc2}, where the same finiteness condition appears in a slightly different guise.) 
An anonymous referee of a precursor to this paper pointed out that with only a little extra work the second of these assumptions can be further weakened to the following: for every $x\in M$ there exist $e_x,f_x\in M$ such that $e_x x=x f_x=x$ and the intervals $[J_x,J_{e_x}]$, $[J_x,J_{f_x}]$ are finite.
This then provides a single proof for Golubov's Theorem \ref{thm10} (ii) and our Theorem \ref{thm16}.
\end{remark}

\section*{Acknowledgements}

The first author was supported by an EPSRC Postdoctoral Fellowship EP/E043194/1 held at the University of St Andrews, Scotland. 

\bibliographystyle{abbrv}

\begin{thebibliography}{10}

\bibitem{AdamsGratzer}
M.E. Adams \and G. Gr\"{a}tzer.
\newblock `Free products of residually finite lattices are residually finite',
\newblock {\em Algebra Universalis} 8 (1978) 262--263. 

\bibitem{Auinger}
K. Auinger. 
\newblock `Residual finiteness of free products of combinatorial strict inverse semigroups',  
\newblock {\em Proc. Roy. Soc. Edinburgh Sect. A}  124  (1994) 137--147. 

\bibitem{BahturinZaicev}
Yu.A. Bahturin \and M.V. Zaicev. 
\newblock `Residual finiteness of color Lie superalgebras', 
\newblock {\em   Trans. Amer. Math. Soc.}  337  (1993) 159--180.

\bibitem{BanaschewskiNelson}
B. Banaschewski \and E. Nelson. 
\newblock `On residual finiteness and finite embeddability', 
\newblock {\em Algebra Universalis}  2  (1972) 361--364.

\bibitem{Carlisle}
W.H. Carlisle. 
\newblock `Residual finiteness of finitely generated commutative semigroups',
\newblock {\em Pacific J. Math. }  36  (1971) 99--101.

\bibitem{CliffordAndPreston}
A.H. Clifford \and G.B. Preston.
\newblock {\em The algebraic theory of semigroups} Vol. II
\newblock (Mathematical Surveys, No. 7., American Mathematical Society,
  Providence, R.I., 1967).

\bibitem{Evans}
T. Evans.  
\newblock `Residually finite semigroups of endomorphisms', 
\newblock {\em J. London Math. Soc.} 2  (1970) 719--721. 

\bibitem{Evans1}
T. Evans. 
\newblock `Some connections between residual finiteness, finite embeddability and the word problem', 
\newblock {\em  J. London Math. Soc.}  1  (1969) 399--403.

\bibitem{Faith}
C. Faith. 
\newblock `Note on residually finite rings',
\newblock {\em Comm. Algebra}  28  (2000) 4223--4226. 

\bibitem{Golubov3}
E.A. Golubov.  
\newblock `Finite separability in semigroups' (Russian),
\newblock {\it Sibirsk. Mat. \v Z. } 11  (1970) 1247--1263. 
(English translation: {\em Siberian Math. J.} 11 (1970) 920--931.)

\bibitem{Golubov4}
E.A. Golubov. 
\newblock `Free product and wreath product of finitely approximable semigroups' (Russian), 
\newblock {\em Ural. Gos. Univ. Mat. Zap.}  8  (1971) 3--15.

\bibitem{Golubov2}
E.A. Golubov.
\newblock `Finitely separated and finitely approximable completely $0$-simple semigroups' (Russian), 
\newblock {\em Mat. Zametki}  12  (1972) 381--392.
(English translation: {\em Math. Notes} 12 (1972)  660--665 (1973).)

\bibitem{Golubov1}
E.A. Golubov.
\newblock `Finitely approximable regular semigroups' (Russian),
\newblock {\em Mat. Zametki}  17 (1975) 423--432.
(English translation: {\em Math. Notes} 17 (1975) 247--251.)

\bibitem{GolubovSapir}
E.A. Golubov \and M.V. Sapir.
\newblock `Ideal extensions and commutative bundles of finitely approximable semigroups' (Russian), 
\newblock {\em Izv. Vyssh. Uchebn. Zaved. Mat.}  7 (1979) 22--30.
(English translation: Soviet Math. (Iz. VUZ) 23 (1979) 21--31.)

\bibitem{GolubovSapir1}
E.A. Golubov \and M.V. Sapir.
\newblock `Residual finiteness of some constructions of universal algebras' (Russian),
\newblock {\em Ural. Gos. Univ. Mat. Zap. }  11  (1978) 26--40.

\bibitem{Gray1}
R.~Gray \and N.~Ru\v{s}kuc.
\newblock `Green index and finiteness conditions for semigroups',
\newblock {\em J. Algebra} 320 (2008) 3145--3164.

\bibitem{Gray2}
R.~Gray \and N.~Ru\v{s}kuc.
\newblock `On residual finiteness of direct products of algebraic systems',
\newblock {\em Monatsh. Math.} 158 (2009), 63--69.

\bibitem{Gruenberg1}
K.W. Gruenberg. 
\newblock `Residual properties of infinite soluble groups', 
\newblock {\em Proc. London Math. Soc.} 7 (1957) 29--62.

\bibitem{PHall1}
P. Hall.
\newblock `On the finiteness of certain soluble groups', 
\newblock {\em Proc. London Math. Soc.} 9 (1959), 595--622.

\bibitem{Hirsch1}
K. Hirsch.
\newblock `On infinite soluble groups. III', 
\newblock {\em Proc. London Math. Soc.} 49 (1946) 184--194.

\bibitem{Howie1}
J.M. Howie.
\newblock {\em Fundamentals of semigroup theory},
\newblock (L.M.S. Monographs 7, Academic Press, Harcourt Brace Jovanovich Publishers, London, 1995.)

\bibitem{Arbib}
K. Krohn, J. Rhodes \and  B. Tilson.
\newblock {\em Algebraic theory of machines, languages, and semigroups},
\newblock Edited by Michael A. Arbib. With a major contribution by Kenneth Krohn and John L. Rhodes, Academic Press, New York, 1968.

\bibitem{Lallement}
G. Lallement. 
\newblock `On nilpotency and residual finiteness in semigroups',  
\newblock {\em Pacific J. Math.}  42  (1972) 693--700. 

\bibitem{Lallement2}
G.~Lallement.
\newblock {\em Semigroups and combinatorial applications},
\newblock (John Wiley \&\ Sons, New York, 1979.)

\bibitem{LallementRosaz}
G. Lallement \and L. Rosaz.
\newblock `Residual finiteness of a class of semigroups presented by a single relation',
\newblock {\em Semigroup Forum}  48  (1994) 169--179.

\bibitem{Magnus1}
W. Magnus. 
\newblock `Residually finite groups',  
\newblock {\em Bull. Amer. Math. Soc.}  75 (1969) 305--316. 

\bibitem{Malcev1}
A. Malcev. 
\newblock `On isomorphic matrix representations of infinite groups' (Russian), 
\newblock {\em Mat. Sb.} 8 (1940) 405--422.

\bibitem{PremetSemenov}
A.A. Premet \and K.N. Semenov. 
\newblock `Varieties of residually finite Lie algebras' (Russian), 
\newblock {\em Mat. Sb. (N.S.)}  137 (1988) 103--113, 144.
(English translation:  Math. USSR-Sb.  65  (1990) 109--118.)

\bibitem{Ruskuc2}
N.~Ru{\v{s}}kuc.
\newblock `On finite presentability of monoids and their Sch\"utzenberger
  groups',
\newblock {\em Pacific J. Math.} 195 (2000) 487--509.

\bibitem{Ruskuc&Thomas}
N.~Ru{\v{s}}kuc \and R.~M. Thomas.
\newblock `Syntactic and Rees indices of subsemigroups',
\newblock {\em J. Algebra} 205 (1998) 435--450.

\bibitem{Varadarajan}
K. Varadarajan. 
\newblock `Rings with all modules residually finite', 
\newblock {\em  Proc. Indian Acad. Sci. Math. Sci.} 109  (1999), 345--351. 

\bibitem{Zaicev}
M.V. Zaicev. 
\newblock `Residual finiteness and the Noethericity of finitely generated Lie algebras' (Russian),
\newblock {\em Mat. Sb. (N.S.)}  136  (1988)  500--509, 591. 
(English translation: Math. USSR-Sb.  64  (1989) 495--504.)
\end{thebibliography}

\end{document}